\documentclass{article}
\usepackage{amssymb,amsfonts,amsmath}

\title{Geometric Methods for Improving the Upper 
Bounds on the Number of Rational Points on Algebraic 
Curves over Finite Fields\\}

\author{Kristin Lauter\\
with an Appendix by J-P.Serre} 
\date{}

\def\F{\mathbb{F}}

\def\Z{\mathbb{Z}}
\def\Q{\mathbb{Q}}
\newcommand{\R}{\mathbb{R}}

\newtheorem{thm}{Theorem}
\newtheorem{lem}{Lemma}

\newtheorem{cor}{Corollary}
\newtheorem{exa}{Example}
\newtheorem{prop}{Proposition}
\newtheorem{rem}{Remark}

\begin{document}
\maketitle


\noindent
\begin{abstract}
Currently, the best upper bounds on the number of rational points
on an absolutely irreducible, smooth, projective algebraic curve  of genus $g$
defined over a finite field $\F_q$ come either from Serre's
refinement of the Weil bound if the genus is small compared to
$q$, or from Oesterl\'e's optimization of the explicit 
formulae method if the genus is large. 

This paper presents three methods for improving these bounds.
The arguments used are the indecomposability of the theta divisor
of a curve, Galois descent, and Honda-Tate theory.  Examples of
improvements on the bounds include lowering them for a wide range of small
genus when $q=2^3, 2^5 ,2^{13}, 3^3, 3^5, 5^3, 5^7$, 
and when $q=2^{2s}$, $s>1$. For
large genera, isolated improvements are obtained for $q=3,8,9$. 

\end{abstract}

\section{Introduction}

This paper presents several methods and results for improving the
upper bounds on the number of rational points on curves over
finite fields.  The first upper bound was discovered in the 1940s by 
Andr\'e Weil as a direct result of proving the 
Riemann Hypothesis for curves.  
Weil showed that the number of rational points, $N$,
on a smooth curve of genus $g$ over the  field $\F_q$ satisfies
the inequality 
\begin{center}
(W) \hspace{20pt} $N \le q+1+2g\sqrt{q}.$
\end{center}
For several decades, number theorists assumed that the Weil
bound was optimal until 1973 when Stark \cite{Stark} improved the bound
by two in a particular case.  

Dramatic improvements began in the 1980s, after Goppa discovered 
that curves over finite fields with many rational points 
could be used to construct efficient error-correcting codes. 
In 1981, Ihara found that equality in the Weil bound can only be
achieved if the genus satisfies
$$g \le \frac{\sqrt{q}(\sqrt{q}-1)}{2}.$$
In \cite{Se}, Serre proved a refined version of the Weil bound:
\begin{center}
(SW) \hspace{20pt} $ N \le q+1+gm, \quad m = [2\sqrt{q}],$
\end{center}
where $[x]$ is the greatest integer part of $x$.
The (SW) bound is the same as (W) if $e$ is even.
For small genus ($g \le \frac{\sqrt{q}(\sqrt{q}-1)}{2}$),
(SW) is the best upper bound known in many cases.
For large genus ($g > \frac{\sqrt{q}(\sqrt{q}-1)}{2}$), 
Serre introduced the so-called explicit formulae method which Oesterl\'e
optimized (see \cite{Se2}). Since its discovery,
the explicit formulae method has provided 
the best upper bounds in large genus. 

Several improvements to (SW)
in the small genus range are already known, notably 
Stark's improvement for $q=13$, $g=2$, (see \cite{Stark});
Serre's generalization for all $g \ge 2$ and $q$ of the form
$q=x^2+1$ or $q=x^2+x+1$, (see \cite{Se2});
Voloch's improvement for $g=3$, $q \le 25$, (see \cite{Se2});
and the improvement when $q$ is a square
and the genus is in the range 
$\frac{(\sqrt{q}-1)^2}{4} < g < \frac{q-\sqrt{q}}{2}$  (see \cite{FT}).

The purpose of this paper is to present three geometric methods for
further improving the bounds. 
The three methods are applications
of Galois descent, Honda-Tate theory, and arguments
on endomorphisms of the Jacobian of a curve.
Examples of improvements obtained via these methods were announced 
in \cite{Lacr}.
For a wide range of small genus we have improved the bounds
by two when $q= 2^3, 2^5, 2^{13}, 3^3, 3^5, 5^3, 5^7$,
and by one  when $q=2^{2s}$, $s>1$.
For large genus, we obtain isolated improvements for $q=3, 2^3, 3^2$.

This paper is organized according to defect.  
A curve has {\it defect $k$} if it fails to meet (SW) by $k$.  
In Section \ref{defect},  we review Serre's derivation of the
list of possible zeta functions for curves of defect $0$, $1$, and $2$.
In Section \ref{descent}, we use Galois descent to treat the
defect $0$ case for $q= 2^3, 2^5, 2^{13}, 3^3, 3^5, 5^3, 5^7$.
In Section \ref{honda}, we use Honda-Tate theory to treat the
defect $2$ case for  $q=2^{2s}$, $s>1$.
Finally, in Section \ref{zeta}, we generate lists of possible
zeta functions for some higher defect cases to improve
several bounds for $q=3,8,9$.  The largest defect we are able to
treat is one case of defect $8$.

The results of this paper explain why 
numerous construction attempts have failed 
to produce curves meeting the explicit formulae bounds in many cases.
Furthermore, the findings presented here
suggest that many more improvements on the bounds can
be made by using a combination of these and other
geometric methods.

\noindent
{\bf Acknowlegements:}  I would like to thank J-P. Serre, James Milne, 
and  Ren\'e Schoof for their generous help and enthusiasm.
In particular, I am grateful to J-P. Serre for suggesting
revisions and for his appendix to improve the proof
of Lemma \ref{deslem}.  Also many thanks to Ren\'e Schoof,
Michael Bennett, and James McLaughlin for pointing out
solutions to the diophantine equations in Section  \ref{solutions}.

\section{Defect $k$}
\label{defect}

\noindent
NOTATION.   Let $q=p^e$, with $p$ prime, $e \ge 1$.
By a curve over $\F_q$, we mean a smooth, projective,
absolutely irreducible curve.  For such a curve, $C$,
let $g=g(C)$ denote the genus, and $N=N(C)$ denote the
number of rational points over $\F_q$.

In this section we investigate the possibilities for 
{\it defect $k$} curves.  

\vspace{10pt} 
\noindent
{\bf Definition } A curve $C$ has {\it defect $k$} if $N(C)=q+1+gm-k$,
$m=[2\sqrt{q}]$. 
\vspace{10pt} 

We explain here the idea used in \cite{Se2} to generate
the list of possible zeta functions for all $(q,g)$
when $N=q+1+gm-k$.
First define the set 
\begin{center}
$F_k = \{ t^d - a_1 t^{d-1} + \dots + a_d \in \Z[t] \mid a_1 = d+k$
and all roots are real $>0$ \}.
\end{center}
Let $F_k^{irred}$ be the subset of $F_k$ consisting of irreducible
polynomials.  Then it follows from Siegel's theorem that
$F_k^{irred}$ is a finite set for $k \ge 0$.
For a fixed $d$, the set of elements  of $F_k$ of
degree $d$ is also finite, and can be listed by taking
all products of elements $f_j$ of degree $d_j$ in $F_{k_j}^{irred}$
such that $\sum k_j =  k$ and $\sum d_j = d$.
In \cite{Smyth},  Smyth produced
complete lists of the sets $F_k^{irred}$, for $k \le 6$.

We say that a curve has {\it zeta function of type $(x_1,\dots,x_g)$}
if $\{\alpha_i, \bar{\alpha_i} \}$ is the family of $g$ conjugate
pairs of eigenvalues of Frobenius acting on the Jacobian of the curve, and  
$x_i = -(\alpha_i + \bar{\alpha_i}), \quad i=1,\dots,g.$
The $m+1-x_i$ are totally positive algebraic integers,
so if $$\sum_{i=1}^{g} x_i = gm-k,$$
then
$$P(t) = \prod_{i=1}^{g}(t-(m+1-x_i)) \in F_k,$$ since
$\deg P =g$, and  $a_1 = g+k$.

The numerator of the zeta function of a curve is determined
by $\{x_i\}$, so a list of possible zeta functions
for curves of defect $k$ and genus $g$ can  be imported
from the lists in \cite{Smyth} for  $k \le 6$.
For $k \le 2$, the lists were contained in \cite{Se2},
and we recall them in Table 1 for convenience.


\begin{table}[h]
\begin{center}

\caption{Possibilities for $(x_1,\dots,x_g)$ for defect $k$, 
with genus restriction}

\begin{tabular}{|c|c|c|}
\hline
k & $(x_1,\dots,x_g)$ & g \\

\hline
$0$    & $(m,\dots,m)$ &  \\
\hline
$1$    & $(m,\dots,m,m-1)$ & $g \ge 1$ \\
& $(m,\dots,m,m+\frac{-1+\sqrt{5}}{2},m+\frac{-1-\sqrt{5}}{2})$ & $g \ge 2$\\
\hline
$2$    & $(m,\dots,m,m-2)$ & $g \ge 1$ \\
       & $(m,\dots,m,m-1,m-1)$ & $g \ge 2$ \\
& $(m,\dots,m,m+\sqrt{2}-1,m-\sqrt{2}-1)$ & $g \ge 2$ \\
& $(m,\dots,m,m+\sqrt{3}-1,m-\sqrt{3}-1)$ & $g \ge 2$ \\
& $(m,\dots,m, m-1,m+\frac{-1+\sqrt{5}}{2},m+\frac{-1-\sqrt{5}}{2})$ & $g \ge 3$\\
& $(m,\dots,m,m+\frac{-1+\sqrt{5}}{2},m+\frac{-1-\sqrt{5}}{2},
m+\frac{-1+\sqrt{5}}{2},m+\frac{-1-\sqrt{5}}{2})$ & $g \ge 4$ \\
& $(m,\dots,m,m+1-4\cos^2\frac{\pi}{7},m+1-4\cos^2\frac{2\pi}{7},
            m+1-4\cos^2\frac{3\pi}{7})$ & $g \ge 3$\\
\hline

\end{tabular}
\end{center}
\end{table}

For $k=3$ the table would have $25$  entries.

For each pair $(q,g)$, an entry in Table 1 might
not correspond to the zeta  function of  a curve for 
a number of possible reasons.
Here are three such reasons from \cite{Se3}.

\vspace{10pt} 
\noindent
{\bf (2.1)}  Since the eigenvalues of Frobenius of a curve
have absolute value $\sqrt{q}$, the $x_i$ must 
satisfy $ |x_i| \le 2\sqrt{q}$.  Any entry in the
table not satifying this condition for all $i$ can be eliminated.
Write $2\sqrt{q} = m + \{2\sqrt{q}\}$, where $\{x\}$ denotes
the fractional part of $x$. Then for example the fourth
entry for defect $2$ is only possible if $\sqrt{3}-1 \le \{2\sqrt{q}\}$.

\vspace{10pt} 
\noindent
{\bf (2.2)} The zeta function of the curve can be expressed in terms of the
$\{x_i\}$ for each entry in the table.  For example
if the Jacobian of the curve is isogenous over  $\F_q$
to a product of elliptic curves,
then $$L(t)=\prod_{i=1}^{g} (1 + x_it+qt^2),$$
and the zeta function is $$Z(t) = \frac{L(t)}{(1-t)(1-qt)}.$$
Re-writing in the form $$1/Z(t) = \prod_{i=1}^{\infty}(1-t^i)^{a_i},$$
we must have $a_i \ge 0$ for all $i$, since $a_i$ is the number of places
of degree $i$ on the curve.
So any entry which does not satisfy this condition can be eliminated.

\vspace{10pt} 
\noindent
{\bf (2.3)} An entry does not correspond to a curve if the Jacobian 
admits a non-trivial decomposition into a product as a polarized
abelian variety.  This condition eliminates any entry for
which the set $\{x_i\}$ can be partitioned into two non-empty subsets,
$I$ and $J$, such that each set is permuted by Gal($\bar{\Q}/\Q$)
and such that the difference between any element of $I$ and any
element of $J$ is a unit.  The proof of this fact can be found in
\cite{Se2} or as Lemma 1 in \cite{La1}.

\vspace{10pt} 
\noindent
A combination of applications of these three reasons provides
the following facts.

\begin{prop} \label{Ihara}
The (SW) bound can only be attained if
$$g \le \frac{q^2-q}{m+m^2-2q}.$$
\end{prop} 
\noindent
Proof: Suppose a curve of genus $g$ attains (SW) so
that $x_i = m$, for all $i$.  The coefficients of the polynomial
$(T+m)^g$ can be computed in two ways: as binomial coefficients
or via Newton's relations between the elementary symmetric
functions, $\{b_n\}$,  and the power functions,
$$s_n = \sum_{i=1}^g (\alpha_i + \bar{\alpha_i})^n.$$
Using the identity 
$$b_2 = \frac{1}{2}(s_1^2 -s_2),$$
and equating the coefficients of the $g-2$ term computed 
in the two ways yields:
$${g \choose 2}m^2 =\frac{1}{2}((gm)^2-(q^2+1-(q+1+gm+2a_2)+2gq)).$$
By reason (2.2), we must have $a_2 \ge 0$, so rearranging yields
the desired inequality.

\begin{rem}
{\rm Note that Proposition 1 generalizes Ihara's result \cite{Ih} that 
the Weil bound cannot be met unless $g \le (q-\sqrt{q})/2.$}
\end{rem}

\begin{prop} \label{defect1}
There are no defect $1$ curves of genus $g > 2$.
\end{prop} 
\noindent
Proof:  This fact was observed in \cite{Se2}, due to reason (2.3), since
both entries for defect $1$ curves can be suitably partitioned
if $g>2$. In the case $g=2$, defect $1$  is only possible if
$\frac{\sqrt{5}-1}{2}  \le \{2\sqrt{q}\}.$

\begin{prop} \label{defect2}
If $e$ is even, then the only defect $2$ curves with genus
$g>2$ have zeta function of type $(m,\dots,m,m-2)$. 
\end{prop} 

\noindent
Proof:  Since $2\sqrt{q} = m$, it follows from reason (2.1)
that the only possibilities are $$(m,\dots,m,m-2)$$ and
$$(m,\dots,m,m-1,m-1).$$  If $g> 2$, then $(m,\dots,m,m-1,m-1)$
is not possible by reason (2.3).

\begin{prop} Let $g \ge 3$, $g \ne 4$.  
If $q$ satisfies $\{2\sqrt{q}\} < \sqrt{3}-1$,
then defect $2$ is only possible if 
$$g \le \frac{q^2-q-2+4m}{m+m^2-2q}.$$
(If $g=4$, then the same conclusion holds if $q$ also
satisfies  $\{2\sqrt{q}\} < \frac{\sqrt{5}-1}{2}).$
\end{prop}
\noindent
Proof: If $g \ge 5$, then by reason (2.3), the only possibilities
are $$(m,\dots,m,m-2)$$ and $$(m,\dots,m,m+\sqrt{3}-1,m-\sqrt{3}-1).$$
The second possibility is eliminated by reason (2.1) since
$\{2\sqrt{q}\} < \sqrt{3}-1$.  The condition on the
genus comes from a computation similar to the one in Proposition
\ref{Ihara}.  Computing the coefficient of the $g-2$ term
in $$(T+m)^{g-1}(T+(m-2))$$ in two different ways and equating
yields:
$$\frac{1}{2}(g-1)(gm^2-4m)=\frac{1}{2}((gm-2)^2-(q^2+1-(q-1+gm+2a_2)+2gq)).$$
Since $a_2 \ge 0$ by reason (2.2), we obtain the stated restriction on the
genus.

If $g=3$,  the last entry in Table 1 is eliminated by reason (2.1) since
$$\{2\sqrt{q}\} < \sqrt{3}-1  < 1-4\cos^2\frac{3\pi}{7}. $$
If  $g=4$, we must also assume that  $\{2\sqrt{q}\} < \frac{\sqrt{5}-1}{2}.$

\section{Galois Descent}
\label{descent}

\vspace{10pt}

When the genus satisfies the inequality of Proposition \ref{Ihara},
we say that the genus is {\it small} (compared to $q$).
For small genus,  (SW) is often the best upper bound known.
Here are some cases where we can improve it by $2$.

\begin{thm} The (SW) bound cannot be met in the following cases:
\label{galois}

\begin{center}
$q=2^3, \quad  4 \le g $, 

$q=2^5, \quad 3 \le g$, 

$q=2^{13}, \quad 4 \le g$, 

$q=3^3, \quad 3 \le g$, 

$q=3^5, \quad 4 \le g$,
 
$q=5^3, \quad 4 \le g$,

$q=5^7, \quad 7 \le g$.

\end{center}
\end{thm}

\begin{rem}
{\rm Note that Serre had already deduced the result in the case $q=27$, 
$g=3$, by using Hermitian modules.  The result for the case $q=243$,
$g=3$  also follows from that argument.
The theorem does not extend to the case $q=8$, $g=3$, since the Klein curve
has $24$ rational points in that case.}  
\end{rem}

\noindent
Proof:  
Theorem 1 is proved by supposing that a curve meeting (SW)
over $\F_q$ {\it does} exist, and using Galois descent to 
produce an $\F_p$-structure for the curve which leads to a contradiction
for the stated cases.
The proof rests on the following descent lemma which was used in
\cite{Se2} to resolve the genus $2$ case.

\begin{lem} \label{deslem}
Let $X$ be a curve over $\F_q$, $q=p^e$, $p$ prime, $e$ odd, 
of genus $g$, $g \ge 2$, with eigenvalues of Frobenius 
$\{\pi, \bar{\pi}\}$ repeated $g$ times.
If $$\pi = \sigma^e, \quad { with } \quad 
\sigma \in \Z[\pi],$$
then $X$ has an $\F_p$-structure with Frobenius endomorphism $\sigma$.
\end{lem}

\noindent Proof of Lemma \ref{deslem}: 
The idea of the proof is as follows.  All details are contained
in the appendix.  If the Jacobian of the curve descends to  $\F_p$
{\it with} its polarization, then the curve also descends by the 
precise version of the Torelli theorem which is stated in the
first section of the appendix.  In order to descend the Jacobian,
it is necessary and sufficient ( \cite{Se4}, Prop. 2, p.110)
that $\sigma$ factors as $$\sigma = \phi \circ \theta,$$ where
$\theta$ is the relative Frobenius map and $\phi$ is a biregular
isomorphism.  The fact that $\sigma$ satisfies this condition
if $\sigma \in \Z[\pi]$ with $\pi = \sigma^e$ and $q=p^e$ with
$p$ prime is shown in Theorem $6'$ of the appendix.  To show that
the principal polarization also descends, it is necessary and
sufficient to show that $$\sigma \sigma ' = p,$$ where $\sigma'$
is the endomorphism obtained from $\sigma$ by applying the involution
associated to the polarization.  The fact that  $\sigma \sigma ' = p$
if  $\sigma \in \Z[\pi]$ is shown in the corollary to Theorem 8
of the appendix.
$\square$

To complete the proof of Theorem 1, we first establish that in the cases
stated in the theorem, a curve meeting (SW) would satisfy 
the hypotheses of Lemma \ref{deslem}.
If $N(X)=q+1+gm$, then the eigenvalues are $\{\pi, \bar{\pi}\}$, 
repeated $g$ times, with $$\pi = \frac{-m \pm \sqrt{m^2-4q}}{2}.$$

\begin{itemize}
\item
$q=2^3$, \quad $m=5$, \quad $\pi = \frac{-5 - \sqrt{-7}}{2} = \sigma^3$, 
\quad $\sigma = \frac{1 + \sqrt{-7}}{2}$, \quad $\sigma = -\pi -2$.

\item
$q=2^5$, \quad $m=11$, \quad $\pi = \frac{-11 + \sqrt{-7}}{2} = \sigma^5$, 
\quad $\sigma = \frac{-1 - \sqrt{-7}}{2}$, \quad $\sigma = -\pi-6 $.

\item
$q=2^{13}$,  $m=181$, \quad $\pi = \frac{-181-\sqrt{-7}}{2}=\sigma^{13}$, 
\quad $\sigma = \frac{1+\sqrt{-7}}{2}$, \quad $\sigma = -\pi-90 $.

\item
$q=3^3$, \quad $m=10$, \quad $\pi = -5+\sqrt{-2} = \sigma^3$, 
\quad $\sigma = 1+\sqrt{-2}$, \quad $\sigma = \pi+6 $.

\item
$q=3^5$, \quad $m=31$, \quad $\pi = \frac{-31 - \sqrt{-11}}{2} = \sigma^5$, 
\quad $\sigma = \frac{-1 - \sqrt{-11}}{2}$, \quad $\sigma = \pi+15 $.

\item
$q=5^3$, \quad $m=22$, \quad $\pi = \frac{-22 - \sqrt{-16}}{2} = \sigma^3$, 
\quad $\sigma = 1+2i$, \quad $\sigma = -\pi-10 $.

\item
$q=5^7$, \quad $m=559$, \quad $\pi = \frac{-559 + \sqrt{-19}}{2} = \sigma^7$, 
$\sigma = \frac{1 + \sqrt{-19}}{2}$, \quad $\sigma = \pi+280 $.

\end{itemize}
Thus $X$ has an $\F_p$-structure with Frobenius $\sigma$, 
so we examine the number of rational
points over $\F_p$ or an extension of $\F_p$,
$$ \#X(\F_{p^e}) = p^e+1-g{\rm Tr}(\sigma^e).$$

\begin{itemize}
\item
$q=2^3$:  Over $\F_2$, $\#X(\F_2) = 2+1-g$, which is impossible
for $g \ge 4$.  In addition, $\#X(\F_4) = 4+1+3g$, which is possible
for $g=3$, but not for $g=2$ or $g \ge 4$.  This gives another
proof of the fact from \cite{Se2} that (SW) cannot be met for
$q=8$, $g=2$.
\item
$q=2^5$: Over $\F_8$, $\#X(\F_8) = 8+1-5g$, which is impossible
for $g \ge 2$.
\item
$q=2^{13}$: Over $\F_2$, $\#X(\F_2) = 2+1-g$, which is impossible
for $g \ge 4$.
\item
$q=3^3$: Over $\F_3$, $\#X(\F_3) = 3+1-2g$, which is impossible
for $g \ge 3$.
\item
$q=3^5$: Over $\F_{27}$, $\#X(\F_{27}) = 27+1-8g$, which is impossible
for $g \ge 4$.
\item
$q=5^3$: Over $\F_5$, $\#X(\F_5) = 5+1-2g$, which is impossible
for $g \ge 4$.
\item
$q=5^7$: Over $\F_5$, $\#X(\F_5) = 5+1-g$, which is impossible
for $g \ge 7$.
\end{itemize}
This completes the proof of Theorem 1. $\square$

For fixed $g$ and $q$, let $N_q(g)$ denote the maximum
of $N(C)$ as $C$ runs through all curves of genus $g$ over $\F_q$.

\begin{cor} If $q= 2^3, 2^{13}, 3^5, 5^3$ (resp. $q= 2^5, 3^3$), and
$g \ge 4$ (resp. $g \ge 3$), then 
$$N_q(g) \le q-1+gm.$$
\end{cor}
Proof: (SW) cannot be met by Theorem 1, 
and Proposition \ref{defect1} implies that defect $1$
is impossible.

\begin{exa}

$$N_8(4) \le 27,$$ $$N_{32}(3) \le 64,$$
$$N_{27}(3) \le 56,$$
$$N_{27}(4) \le 66.$$

\end{exa}

\subsection{A pair of Diophantine equations}
\label{solutions}

Cases where Lemma \ref{deslem} improves the upper bounds for the
number of rational points on curves over a finite field $\F_q$
correspond to integer solutions to a pair of diophantine equations:
$${\bf (3.1)}  \quad \quad \quad    x^2+d=4p,$$
$${\bf (3.2)}  \quad \quad \quad y^2+d=4p^e,$$
where $d$ is positive, $e$ is odd, $p$ is prime, and $q=p^e$.
Provided that $d<2y+1$, we have $y=m=[\sqrt{4q}]$,
and so $-d=m^2-4q$ as in the instances of Theorem \ref{galois}.
The correspondence is expressed by the following lemma.

\begin{lem} \label{nagell}
A solution  $(x,y,d,p,e)$ to the pair of equations
(3.1) and (3.2) with $3 < d < 2y+1$, $d$ 
square-free, corresponds to a pair of algebraic
integers $\pi$ and $\sigma$ which satisfy the
conditions of Lemma \ref{deslem}: 
$$\pi = \frac{-y \pm \sqrt{-d}}{2}, \quad
\sigma =  \frac{x \pm \sqrt{-d}}{2}, \quad
\sigma^e=\pm \pi, \quad \sigma \in \Z[\pi], \quad 
\pi\bar{\pi}=q, \quad \sigma\bar{\sigma}=p.$$
\end{lem}

\noindent
Proof: 
Given  a solution to (3.1) and (3.2), we must show that
$\pi$ and $\sigma$ as defined in the statement satisfy
the required properties.  All properties are immediate except
$\sigma^e=\pm \pi$. This follows from the fact that $R$,
the ring of integers in  $\Q(\pi)$, is a unique factorization
domain  with only the trivial  units $\{\pm 1\}$.
In fact, $p$ splits in $R$ as $p=\sigma\bar{\sigma}$,
but $p$ does not divide $\pi$.  Since $p^e=\pi\bar{\pi}$
we  must have $\pi$ or $\bar{\pi}$ associated to $\sigma^e$.
$\square$

\begin{exa} {\rm 
Note that if $q=7^3$, we have the following solution
to (3.1) and (3.2): $x=5$, $y=37$, $d=3$.
This solution does not correspond to a  pair satisfying
the conditions of Lemma \ref{deslem} however,
since there are non-trivial units in the ring of integers
of $\Q(\sqrt{-3})$. Let $u= \frac{1 + \sqrt{-3}}{2}$. Then
if $\pi =  \frac{-37 - \sqrt{-3}}{2}$ and
$\sigma =  \frac{-5 + \sqrt{-3}}{2}$, we have 
$\pi = u \cdot \sigma^3$. }
\end{exa}

In \cite{Se2}, Serre deduced that (SW) cannot be met
when $q$ is of the form  $q=x^2 + 1$ or $x^2 + x + 1$
and the genus is at least $2$. 
Theorem \ref{galois} was discovered as a result of
trying to extend this theorem to prime powers
of the form  $q=x^2 + x+2$.
For $q$ of the form  $q=x^2 + x+2$, we must have $p=2$, and there are only
$5$ such $q$.  They correspond to the famous solutions to
(3.1) and (3.2) when $d=7$ and $e=1,2,3,5,13$, referred
to in this case as the Ramanujan-Nagell equations.

\vspace{10pt} 
In addition to the solutions to (3.1) and (3.2) listed in
Theorem \ref{galois}, there is a family of solutions when $e=3$
which was pointed out by Ren\'e Schoof and Michael
Bennett.  For any integer $k$
such that $p=k^2+1$ is prime, we have a solution
of the form $$x=k, \quad y=k(2k^2+3), \quad
d=3k^2+4, \quad e=3.$$
For $k=1$ and $k=2$, the solutions correspond to
the first and the second-to-last instances of Theorem \ref{galois}.
This (conjecturally) infinite family of solutions 
improves the upper bounds for the number of rational
points on curves for each of the corresponding fields,
due to Lemma \ref{nagell}.

\section{Honda-Tate Theory}
\label{honda}

Since the Weil and Serre bounds coincide
when $q$ is a square, we can consider the defect from the Weil bound.  
For $q=2^{2s}, \quad s>1$ and $g$ small in a certain range,
we can improve the bounds due to the following theorem.

\begin{thm} If $q=2^{2s}$, $s>1$, and $g>2$, 
then there are no defect $2$ curves.
\end{thm}
Proof: The proof of Theorem 2 relies on Honda-Tate theory.
By Proposition \ref{defect2},
for $q$ a square and $g>2$,  the only possibility
for a defect $2$ curve is one with its Jacobian
isogenous to the product of elliptic curves:
$$E_m \times \dots\times E_m \times E_{m-2},$$
where $E_m$ is an elliptic curve with Tr(Frobenius) $=-m$.
By Honda-Tate theory, when $q=2^{2s}$, the only possible
values for the trace of an elliptic curve which are
divisible by the characteristic are (see \cite{Water}, p.536)
$$\{0, \pm \sqrt{q}, \pm 2\sqrt{q} \}.$$
If $s > 1$, then $m-2 = 2\sqrt{q} - 2$ is not on this list,
so such an abelian variety is impossible.

\begin{cor} If  $q=2^{2s}$, $s>1$, and
$\frac{(\sqrt{q}-1)^2}{4} < g < \frac{q-\sqrt{q}}{2},$
then $$N_q(g) \le q-2+gm.$$
\end{cor}
Proof:  Due to a result of Fuhrmann and Torres \cite{FT}
when $q$ is a square, there are no defect $0$ curves
for any $g$ in the interval
$$\frac{(\sqrt{q}-1)^2}{4} < g < \frac{q-\sqrt{q}}{2}.$$
By Proposition \ref{defect1}, there are no defect $1$ curves
for $g > 2$.  By Theorem 2, there are no defect $2$ curves
for $g > 2$.

\begin{exa}
{\rm Theorem 2 leads to the following improvements on the bounds:
$$N_{16}(4) \le 46,$$
$$N_{16}(5) \le 54,$$
$$N_{64}(g) \le 62+16g, \quad {\rm  for } \quad 13 \le g \le 27.$$
} 
\end{exa}

\section{Zeta Functions}
\label{zeta}

When the genus is large, the explicit formulae bounds
force the maximal curves to have defect
bigger than $2$.  In this case, we proceed by using
reasons \#1,2,3 from Section \ref{defect} directly 
to generate lists of possible zeta functions.   
In the following theorem we give several cases where the explicit
formulae bounds can be improved by showing that the lists are empty.

\begin{thm}
The optimal form of the explicit formulae
bounds cannot be met in the following cases:

\begin{center}
$q=3, \quad g=5, \quad N=14$ (defect 5)

$q=3, \quad g=7, \quad N=17$ (defect 8)

$q=9, \quad g=5, \quad N=36$ (defect 4)

$q=8, \quad g=6, \quad N=36$. (defect 3)
\end{center}
\end{thm}
Proof: 
The first two cases were proved in \cite{La1} and \cite{La2}
respectively.  

For $q=9$, $g=5$, the Weil and Serre bounds give $N \le 40$.
The optimal form of the explicit formulae bounds gives
$N \le 36$, which is defect $4$.  Applying the algorithm
from \cite{La2}, we find that the only possibility is
$$(m,m-1,m-1,m-1,m-1),$$ which is impossible by reason \#3.

For $q=8$, $g=6$, the  optimal form of the explicit formulae bounds 
again gives $N \le 36$, which is defect $3$.  The possibilities
are $$(m,m,m,m,m-1,m-2),$$
$$(m,m,m,m-1,m-1,m-1),$$
$$(m,m,m-1,m-1,m-\frac{1+\sqrt{5}}{2},m-\frac{1-\sqrt{5}}{2}),$$
all three of which are impossible by reason \#3.

\newtheorem{thmf}{Th\'eor\`eme}

\bigskip\bigskip

\noindent\centerline{{\large \bf{Appendice}}}
\smallskip

\noindent{ch\`ere Kristin,}                  

\noindent
A propos de la descente du corps de base pour les courbes et
leurs jacobiennes:

\bigskip\noindent
\centerline{\bf 1. Le th\'eor\`eme de Torelli}
\smallskip

Soit $k$ un corps. Par une ``courbe'' sur $k$ j'entends une
courbe projective, lisse, absolument irr\'eductible.
Si $X$ est une telle courbe, son genre $g(X)$ sera not\'e $g$.
On suppose $g>1$.  On note Jac $X$ la jacobienne de $X$ munie
de sa polarisation naturelle $a$, qui est de degr\'e $1$.
Si $X'$ est une autre courbe sur $k$, tout isomorphisme
$f: X \rightarrow X'$ d\'efinit par transport de structure
un isomorphisme $f_J:(J,a) \rightarrow (J',a')$, o\`u $(J',a')$
est la jacobienne de $X'$. Le {\bf th\'eor\`eme de Torelli}
dit que l'on obtient ainsi ``presque'' tous les isomorphismes \\
$(J,a) \rightarrow (J',a')$.  De fa\c con plus pr\'ecise:

\begin{thmf}
Supposons $X$ hyperelliptique. Pour tout isomorphisme
de vari\'et\'es ab\'eliennes polaris\'ees
$$F: (J,a) \rightarrow (J',a'),$$
il existe un isomorphisme $f:X \rightarrow X'$ et un seul tel que $F=f_J$.
\end{thmf}

\begin{thmf}  
Supposons $X$ non hyperelliptique.  Alors, pour tout
isomorphisme $F: (J,a) \rightarrow (J',a'),$ il existe un
isomorphisme $f:X \rightarrow X'$ et un entier $e$ \'egal \`a
$\pm 1$ tel que $F=e \cdot f_J$.  De plus, le couple $(f,e)$
est d\'etermin\'e par $F$ de fa\c con unique.
\end{thmf}

(Noter que $X$ est hyperelliptique si et seulement si il existe un
automorphisme $s$ de $X$ tel que $s_J=-1$.)

\noindent
Les ths.1 et 2 constituent ce que j'ai envie d'appeler la
``forme pr\'ecise'' du th\'eor\`eme de Torelli (la forme impr\'ecise
consistant \`a dire seulement que $X$ et $X'$ sont isomorphes). 
Je ne crois pas que la ``forme pr\'ecise'' se trouve explicitement
dans la litt\'erature.  Toutefois:

---Lorsque $k$ est alg\'ebriquement clos, c'est essentiellement
l'\'enonc\'e d\'emontr\'e par Weil ({\it Oe.} II, [1957a]),
\`a cela pr\`es que Weil choisit un plongement
de $X$ dans sa jacobienne, ce qui introduit des translations
qui n'ont rien \`a voir avec la question.  La d\'emonstration
du th\'eor\`eme de Torelli due \`a Andreotti (et reproduite par exemple
dans Albarello-Cornalba-Griffiths-Harris, Grundlehren 267) ne
donne que la forme impr\'ecise.

---Le cas d'un corps parfait r\'esulte du cas alg\'ebriquement clos
par descente galoisienne standard (gr\^ace \`a {\it l'unicit\'e}
de $f$ ou $(f,e)$).  Le cas d'un corps imparfait r\'esulte de
celui d'un corps parfait : en effet, si $k_1$ est une extension
radicielle de $k$, tout isomorphisme de $X/k_1$ sur $X'/k_1$
est ``d\'efini sur $k$'', i.e. provient d'un isomorphisme
de $X$ sur $X'$ (utiliser le fait que le sch\'ema Isom($X$,$X'$)
est \'etale). D'ailleurs, dans la suite, le cas d'un corps parfait nous 
suffira.

\bigskip\noindent
\centerline{\bf 2. Un corollaire du th\'eor\`eme de Torelli}
\smallskip

C'est l'\'enonc\'e suivant, qui r\'esulte imm\'ediatement des 
ths.1 et 2:

\begin{thmf} 
On a
$$ {\rm Aut}(J,a)
=
\left\{
\begin{array}{cl}
{\rm Aut} X  & \mbox{si $X$ est hyperelliptique} \\
\{\pm 1 \} \times {\rm Aut} X & \mbox{si $X$ n'est pas hyperelliptique. } \\
\end{array}
\right.
$$
\end{thmf}

{\bf Corollaire} {\it Supposons que le groupe fini ${\rm Aut} (J,a)$ 
contienne un \'el\'ement $s$ tel que $s^n=-1$, avec $n$ pair.  Alors $X$ est
hyperelliptique.}

En effet l'existence d'un tel $s$ est incompatible avec la d\'ecomposition \\
$ {\rm Aut}(J,a) =  \{\pm 1 \} \times {\rm Aut} X$.

Cet \'enonc\'e peut \^etre utile pour montrer que certaines courbes, 
construites par la m\'ethode des modules hermitiens, sont
hyperelliptiques.

\bigskip\noindent
\centerline{\bf 3. Descente du corps de base}
\smallskip

On se donne une extension galoisienne finie $k_1/k$, de groupe
de Galois $G$.  Pour \'eviter des indices trop abondants, on
note $X_1$, $J_1$,... une courbe sur $k_1$, sa jacobienne, etc.
On se donne {\it une $k$-structure sur $J_1$, compatible avec
la polarisation}; cela revient \`a se donner
un couple $(J,a)$, o\`u $J$ est une vari\'et\'e ab\'elienne
sur $k$, munie d'une polarisation $a$ d\'efinie sur $k$, et \`a se donner
un isomorphisme de $(J_1,a_1)$ avec $(J,a)/k_1$.  On veut passer
de la jacobienne \`a la courbe.
\begin{thmf} 
Supposons $X_1$ hyperelliptique.  Il existe alors une
$k$-structure unique sur $X_1$, compatible avec sa $k_1$-structure, et
dont la jacobienne est $(J,a)$.
\end{thmf}

Cela r\'esulte du th.1, par descente \`a la Weil.
De fa\c con plus pr\'ecise, si $s$ est un \'el\'ement donn\'e
de $G$, la $k$-structure $(J,a)$ donne un isomorphisme de
$(J_1,a_1)$ sur son $s$-transform\'e $(J_1,a_1)^s$; d'o\`u
par le th.1 un isomorphisme $f_s : X_1 \rightarrow (X_1)^s$.
Ces isomorphismes satisfont \`a la condition de cocycle usuelle.  D'o\`u
la structure cherch\'ee.

Le cas non hyperelliptique est analogue, mais plus amusant:
\begin{thmf}
Supposons $X_1$ non hyperelliptique.  Il existe alors une
$k$-structure $X$ sur $X_1$, compatible avec sa $k_1$-structure, et
un homomorphisme $\epsilon : G \rightarrow \{\pm 1\}$, tel que la
jacobienne de $X$ soit isomorphe \`a la $\epsilon$-tordue 
$(J,a)_{\epsilon}$ de $(J,a)$.
\end{thmf}

\indent
(Par `` $\epsilon$-tordue '' j'entends la vari\'et\'e d\'eduite de
$(J,a)$ par torsion galoisienne relativement \`a 
$\epsilon : G \rightarrow \{\pm 1\} \subset {\rm Aut}(J,a)$.)

\indent
La d\'emonstration est la m\^eme que celle du th 4.
Pour chaque $s \in G$ on a  (cf. th.2) un isomorphisme
$f_s : X_1 \rightarrow (X_1)^s$ ainsi qu'un signe $e_s=\pm 1$.
On d\'efinit alors $\epsilon$ par $\epsilon (s) = e_s$.

{\it Remarque}.  On pourrait s\^urement d\'eduire les ths.
4 et 5 d'un \'enonc\'e portant sur le morphisme de ``champs'': \\
``champ de courbes'' $\rightarrow$ ``champ de vari\'et\'es
ab\'eliennes \`a polarisation principale''. 

\bigskip\noindent
\centerline{\bf 4. Corps finis :  vari\'et\'es ab\'eliennes}
\smallskip

On va s'int\'eresser maintenant au cas o\`u $k$ est un corps fini
\`a $q$ \'el\'ements et $k_1$ une extension finie \`a $q_1$
\'el\'ements, avec $q_1=q^r$, $r>1$. On se donne une courbe $X_1$
sur $k_1$, et l'on d\'esire ``descendre'' son corps de d\'efinition
\`a $k$, comme ci-dessus.  Cela va se faire en trois \'etapes:
\begin{list}{}{}
\item
descente pour les  vari\'et\'es ab\'eliennes;
\item
descente pour les  vari\'et\'es ab\'eliennes polaris\'ees;
\item
descente pour les courbes. 
\end{list}

Occupons-nous du premier cas, i.e. de celui des vari\'et\'es ab\'eliennes.
On se donne une vari\'et\'e ab\'elienne $A_1$ sur $k_1$. Notons $\pi_1$
son endomorphisme de Frobenius.  Une $k$-structure sur $A_1$ est
d\'efinie par son endomorphisme de Frobenius \\ $\pi \in {\rm End} (A_1)$.
Les conditions que $\pi$ doit satisfaire sont les suivantes:
\begin{thmf} 
Pour que $\pi$ d\'efinisse sur $A_1$ une $k$-structure
compatible avec sa $k_1$-structure, il faut et il suffit que:

a) $\pi_1 = \pi^r$, o\`u $r=[k_1:k]$;

b) $\pi$ est nul sur le noyau $N_q$ de l'homomorphisme de Frobenius
absolu \\
$$F_q : A_1 \rightarrow (A_1)^{(q)}.$$ 
\end{thmf}

(Une fa\c con \'equivalente de formuler $b)$ est de dire que, pour
toute fonction rationnelle $h$ sur $A_1$, la fonction $h \circ \pi$ est
la puissance $q$-i\`eme d'une fonction rationnelle.)

La n\'ecessit\'e de ces conditions est imm\'ediate.  La suffisance
r\'esulte par exemple de la prop.2 de 
[8], Chap.VI $\S$ 1, p.110: d'apr\`es cette proposition,
on doit v\'erifier que $\pi$ est de la forme $\sigma \circ F_q$, o\`u
$\sigma$ est un isomorphisme de $(A_1)^{(q)}$ sur $A_1$.  Or $b)$
entra\^{\i}ne que $\pi$ se factorise en $\sigma \circ F_q$, o\`u
$\sigma$ est un homomorphisme de $(A_1)^{(q)}$ dans $A_1$.
Comme les degr\'es de $\pi$ et de $F_q$ sont tous deux \'egaux \`a 
$q^{\text{dim}(A)}$, on voit que $\sigma$ est de degr\'e $1$,
i.e. que c'est un isomorphisme.

{\it Remarque}. On peut donner des exemples o\`u $a)$ est v\'erifi\'ee,
mais pas $b)$. \\  Toutefois:

{\bf Th\'eor\`eme 6$'$} {\it Supposons que la condition $a)$ du th.6 
soit satisfaite.  Faisons les hypoth\`eses suivantes :

c) $q$ est \'egal \`a la caract\'eristique $p$ du corps $k$;

d) il existe un polyn\^ome $P(X)$ \`a coefficients entiers tel que
$\pi$ soit \'egal \`a $P(\pi_1)$.

Alors la condition b) du th.6 est satisfaite.}

Ecrivons $\pi$ sous la forme $a_0+a_1\pi_1+...+a_n\pi_1^n$, avec
$a_i \in \Z$.  L'application tangente \`a $\pi_1$ est nulle.
Il en r\'esulte que l'application
tangente \`a $\pi$ est l'homoth\'etie de rapport $a_0$.
D'apr\`es $a)$ la puissance $r$-i\`eme de cette application est $0$.
Il en r\'esulte que $a_0$ est divisible par $p$, d'o\`u le fait
que l'application tangente \`a $\pi$ est $0$.  Or cela signifie
que $\pi$ s'annule sur $N_p$.  La condition $b)$ est donc satisfaite.

\bigskip\noindent
\centerline{\bf 5. Corps finis :  vari\'et\'es ab\'eliennes polaris\'ees}
\smallskip

On conserve les hypoth\`eses du \S 4, et l'on suppose en outre que
$A_1$ est munie d'une polarisation $a_1$.  On se donne 
$\pi \in \text{End}(A_1)$ satisfaisant aux conditions du th.6,
donc d\'efinissant sur $A_1$ une $k$-structure.  Soit $A$
la vari\'ete ab\'elienne ainsi obtenue.  On d\'esire donner des conditions
permettant d'affirmer que $a_1$ est d\'efinie sur $k$, i.e. provient
d'une polarisation $a$ de $A$.

Je rappelle qu'une polarisation d'une vari\'et\'e ab\'elienne
d\'efinit une involution de l'alg\`ebre $R_{\Q}={\Q} \otimes R$, o\`u
$R= \text{End} (A)$ (cette involution laisse stable $R$ lorsque la
polarisation est degr\'e $1$, ce qui est le cas qui nous int\'eresse
le plus).  Je noterai $x \mapsto x'$ l'involution d\'efinie par
la polarisation $a_1$.  En particulier, $\pi '$ est  d\'efini.
\begin{thmf} 
Pour que la polarisation $a_1$ soit rationnelle sur $k$
(pour la $k$-structure  d\'efinie par $\pi$), il faut et il suffit que
l'on ait

$e)$ $\pi \pi ' =q$.
\end{thmf}
Notons $V$ le ``Verschiebung'' de $A$, i.e. l'unique endomorphisme de
$A$ tel que $\pi V = q$.  La condition $e)$ ci-dessus \'equivaut \`a:

$e')$ $\pi ' = V$.

(Rappel de notations: si $C$ est une vari\'et\'e ab\'elienne,
je note $C^*$ sa duale; de m\^eme, si $h: B \rightarrow C$ 
est un homomorphisme,
je note $h^*$ l'homomorphisme correspondant (``transpos\'e'',
``adjoint'',...) de $C^*$ dans $B^*$.  La polarisation 
$a_1 : A \rightarrow A^*$ est hermitienne: on a $a_1^*=a_1$.
L'involution associ\'ee $x \mapsto x'$ de ${\Q}  \otimes R$
est caract\'eris\'ee par la formule $a_1.x' = x^*.a_1$.)

Il est bien connu que l'endomorphisme de Frobenius de $A^*$
est \'egal \`a $V^*$.  Le morphisme $a_1: A \rightarrow A^*$
est rationnel sur $K$ si et seulement si il commute au Frobenius,
i.e. si et seulement si on a $a_1 \pi = V^* a_1$. En comparant
\`a l'\'equation $a_1 V' = V^*a_1$, on voit que cela revient \`a
$\pi = V'$, i.e. \`a $\pi' = V$.  D'o\`u le th\'eor\`eme.

Les conditions $e)$ et $e')$ peuvent \^etre remplac\'ees par
une condition plus simple:
\begin{thmf} 
La condition $e)$ \'equivaut \`a :

$e'')$ $\pi$ et $\pi'$ commutent.
\end{thmf}

Il est clair que $e) \Rightarrow e'')$.  Pour prouver la r\'eciproque,
il est commode d'utiliser l'alg\`ebre $S = {\bf{R}} \otimes \text{End}(A)$,
et la sous-alg\`ebre $T$ de $S$ engendr\'ee par $\pi$ et $\pi'$.
Vu $e'')$, cette alg\`ebre est commutative, et stable par l'involution
$x \mapsto x'$.  De plus, si l'on pose $v = \pi/q^{1/2}$, et $z=vv'$,
on a  $z^r = \pi^r {\pi'}^r/q^r = (\pi_1 \pi_1')/q^r = 1$.
Or l'alg\`ebre $S$ (et donc aussi l'alg\`ebre $T$) peut \^etre 
munie (voir Mumford) d'une forme lin\'eaire r\'eelle $t$ telle que
$t(yy') > 0$ pour tout $y \ne 0$.  Comme $T$ est commutative,
il en r\'esulte que $T$ se d\'ecompose en produit de corps
isomorphes \`a $\R$ ou ${\bf{C}}$, l'involution \'etant la 
conjugaison complexe.  La formule $z=vv'$ montre que, dans chacun de 
ces corps, $z$ est r\'eel $ > 0$.  Comme d'autre part c'est une racine 
de l'unit\'e, on a $z=1$, ce qui \'equivaut \`a $\pi\pi' = q$.

{\bf Corollaire} {\it La condition d) du th.6$'$ entra\^{\i}ne la 
condition e) du th.7.}

En effet, si $\pi$ est un polyn\^ome en $\pi_1$, $\pi'$ est un polyn\^ome
en $\pi_1'$. Or $\pi_1$ et $\pi_1'$ commutent (puisque $a_1$ est d\'efinie 
sur $k_1$).  Donc $\pi$ et $\pi'$ commutent, et l'on peut appliquer le th.8.

\bigskip
\noindent
\centerline{\bf 6. Corps finis: courbes}
\smallskip

On se donne une courbe $X_1$ sur $k_1$, et l'on note $\pi_1$
l'endomorphisme de Frobenius de sa jacobienne $J_1$.  On se donne
$\pi \in \text{End} (J_1)$, avec $\pi^r = \pi_1$ et l'on cherche
\`a mettre sur $X_1$ une $k$-structure telle que l'endomorphisme
de Frobenius  corres-pondant soit $\pi$.  On suppose que $\pi$
satisfait aux conditions $b)$ et $e)$ des ths. 6 et 7.  Alors:
\begin{thmf}
Si $X_1$ est hyperelliptique, ou si $r$ est impair,
il existe sur $X_1$ une $k$-structure dont le Frobenius est $\pi$.

Sinon, il existe un signe $\epsilon = \pm 1$ et une $k$-structure
sur $X_1$ dont le Frobenius est $\epsilon \pi$.
\end{thmf} 

Vu les ths.6 et 7, il existe une $k$-structure sur $J_1$,
compatible avec sa $k_1$-structure et sa polarisation, pour laquelle le
Frobenius est $\pi$.  Si $X_1$ est hyperelliptique, le th.4 donne
l'existence de la $k$-structure cherch\'ee sur $X_1$.
 Si $X_1$ n'est pas hyperelliptique, le th.5 donne le m\^eme
r\'esultat, \`a cela pr\`es que la $k$-structure de $J_1$ doit
\^etre tordue par un caract\`ere quadratique de Gal($k_1/k$),
qui est un groupe cyclique d'ordre $r$.  Si $r$ est impair, un tel
caract\`ere est trivial; aucune torsion n'est donc n\'ecessaire.  
Si $r$ est pair, il se peut que ce caract\`ere soit l'unique
caract\`ere non trivial; or l'effet d'une telle torsion est de
remplacer $\pi$ par son oppos\'e.  D'\`ou le r\'esultat cherch\'e.

{\it Remarque}.  Dans le cas qui vous int\'eresse, on a $q=p$,
et $\pi$ s'\'ecrit comme polyn\^ome en $\pi_1$ \`a coefficients
entiers.  Les conditions $b)$ et $c)$ sont alors satisfaites, et le 
th.9 s'applique.

C'est ce que l'on voulait.

\bigskip
\noindent{Bien \`a vous}

\bigskip
\noindent{J.-P. Serre}



\end{document}